\documentclass[12pt]{amsart}
\usepackage{amsmath}
\usepackage{amsxtra}
\usepackage{amscd}
\usepackage{amsthm}
\usepackage{amsfonts}
\usepackage{amssymb}
\usepackage{eucal}
\usepackage{graphicx}

\input xy
\xyoption{all}

\theoremstyle{remark}

\theoremstyle{definition}

\numberwithin{thm}{section}
\numberwithin{equation}{section}

\newcommand{\nc}{\newcommand}
\nc{\on}{\operatorname}
\nc{\ch}{\mbox{ch}}
\nc{\Z}{{\mathbb Z}}
\nc{\C}{{\mathbb C}}
\nc{\pone}{{\mathbb P}^1}
\nc{\pa}{\partial}
\nc{\F}{{\mathcal F}}
\nc{\arr}{\rightarrow}
\nc{\larr}{\longrightarrow}
\nc{\al}{\alpha}
\nc{\ri}{\rangle}
\nc{\lef}{\langle}
\nc{\W}{{\mathbb W}}
\nc{\la}{\lambda}
\nc{\ep}{\epsilon}
\nc{\su}{\widehat{{\mathfrak s}{\mathfrak l}}_2}
\nc{\sw}{{\mathfrak s}{\mathfrak l}}
\nc{\g}{{\mathfrak g}}
\nc{\h}{{\mathfrak h}}
\nc{\n}{{\mathfrak n}}
\nc{\N}{\widehat{\n}}
\nc{\G}{\widehat{\g}}
\nc{\De}{\Delta}
\nc{\gt}{\widetilde{\g}}
\nc{\Ga}{\Gamma}
\nc{\one}{{\mathbf 1}}
\nc{\z}{{\mathfrak Z}}
\nc{\La}{\Lambda}
\nc{\wt}{\widetilde}
\nc{\wh}{\widehat}
\nc{\cri}{_{\kappa_c}}
\nc{\kk}{\underline{\mathbf C}}
\nc{\sun}{\widehat{\sw}_N}
\nc{\si}{\sigma}
\nc{\el}{\ell}
\nc{\bi}{\bibitem}
\nc{\om}{\omega}
\nc{\ol}{\overline}
\nc{\ds}{\displaystyle}
\nc{\dzz}{\frac{dz}{z}}
\nc{\Res}{\on{Res}}
\nc{\mc}{\mathcal}
\nc{\Cal}{\mathcal}
\nc{\bb}{{\mathfrak b}}
\nc{\ot}{\otimes}
\nc{\R}{{\mc R}}
\nc{\yy}{{\mc Y}}
\nc{\ga}{\gamma}

\nc{\us}{\underset}
\nc{\opl}{\oplus}

\nc{\Fq}{{\mathbb F}_q}
\nc{\Mq}{{\mathcal M}}
\nc{\Rep}{\on{Rep}}
\nc{\sssec}{\subsubsection}
\nc{\ssec}{\subsection}
\nc{\lan}{\langle}
\nc{\ran}{\rangle}

\nc{\D}{\mathcal D}
\nc{\Vect}{\on{Vect}}
\nc{\ghat}{\G}
\nc{\T}{\mc T}
\nc{\Tloc}{\T^\g_{\on{loc}}}
\nc{\vac}{|0\ran}
\nc{\Wick}{{\mb :}}
\nc{\mb}{\mathbf}
\nc{\delz}{\partial_z}
\nc{\K}{{\mathbb K}}
\nc{\cali}{\mathcal}
\nc{\li}{\mathfrak l}
\nc{\lt}{\widetilde{\li}}
\nc{\astar}{a^*}
\nc{\cA}{{\mc A}}
\nc{\ka}{\kappa}

\nc{\OO}{{\mc O}}
\nc{\AutO}{\on{Aut}\OO}
\nc{\DerO}{\on{Der}\OO}
\nc{\DerpO}{\on{Der}_+\OO}
\nc{\Au}{{\mc A}ut}
\nc{\mf}{\mathfrak}
\nc{\V}{{\mathcal V}}
\nc{\hh}{\wh{\h}}

\nc{\pp}{{\mathfrak p}}
\nc{\mm}{{\mathfrak m}}
\nc{\rr}{{\mathfrak r}}
\nc{\ket}{\rangle}
\nc{\zz}{{\mathfrak z}}
\nc{\gr}{\on{gr}}
\nc{\Spe}{\on{Spec}}
\nc{\rv}{\crho}
\nc{\can}{\on{can}}
\nc{\CC}{\on{Op}_G(D))}
\nc{\Op}{\on{Op}_G(D)}
\nc{\MOp}{\on{MOp}_G(D)}
\nc{\Db}{{\mathbb D}}
\nc{\ww}{w}

\nc{\af}{{\mathbb A}^1}
\nc{\bs}{\backslash}
\nc{\laa}{(\la_i)}
\nc{\zn}{(z_i)}

\nc{\cla}{\check{\la}}
\nc{\cmu}{\check{\mu}}
\nc{\crho}{\check{\rho}}
\nc{\chal}{\check{\al}}
\nc{\cc}{{\mathfrak c}}

\nc{\MM}{{\mathbb M}}

\nc{\ZZ}{{\mc Z}}

\nc{\UU}{{\mathbb U}}

\nc{\Conn}{\on{Conn}(\Omega^{\crho})}
\nc{\Con}{\on{Conn}(\Omega^{-\rho})}
\nc{\Co}{\on{Conn}(\Omega^{\rho})}

\nc{\ppart}{(\!(t)\!)}
\nc{\pparl}{(\!(\la)\!)}
\nc{\zpart}{(\!(z)\!)}
\nc{\ppzi}{(\!(t-z_i)\!)}
\nc{\ppinf}{(\!(t^{-1})\!)}
\nc{\Ind}{\on{Ind}}
\nc{\I}{{\mathbb I}}

\nc{\Bun}{\on{Bun}}

\nc{\gtil}{\wt{\g}}
\nc{\ntil}{\wt{\n}}
\nc{\htil}{\wt{\h}}
\nc{\gbar}{\ol{\g}}
\nc{\nbar}{\ol{\n}}
\nc{\bbar}{\ol{\bb}}
\nc{\lhat}{\wh{\mf l}}

\nc{\ovc}{\overset{\circ}}
\nc{\Gr}{\on{Gr}}

\nc{\AD}{{\mathbb A}}
\nc{\Gm}{{\mathbb G}_m}
\nc{\Ql}{{\mathbb Q}_\ell}

\nc{\Loc}{\on{Loc}}

\nc\GG{\mathbb G}
\nc\Xb{\mathbf X}
\nc\inv{{\rm inv}}
\nc\isom{=}
\nc\Hc{{\mathcal H}}
\nc\ovl{\overline}

\nc\Lc{{\mathcal L}}
\nc\Ec{{\mathcal E}}
\nc\wF{{\mc K}}
\nc\pr{{\rm pr}}
\nc\SL{{\rm SL}}
\nc\PGL{{\rm PGL}}
\nc\GL{{\rm GL}}
\nc\LG{{}^L G}
\nc\Pc{{\mc P}}
\nc\Ac{{\mc A}}
\nc\dv{{\rm div}}
\nc\Fc{{\mc F}}
\nc\M{{\mc M}}

\nc{\und}{\underline}

\nc{\zb}{\ol{z}}
\nc{\wb}{\ol{w}}
\nc{\p}{{\partial}}
\nc{\pf}{\int\hspace*{-3.5mm}\bs}

\begin{document}

\title{Mathematics, Love, and Tattoos}\thanks{Contribution to the
  Proceedings of the Symposium ``Matematica e Cultura'' 2012, Venice
  (ed. M. Emmer) to be published by Springer Verlag. \\ Parts of the
  article are borrowed from my book {\em Love and Math}, which will be
  published by Basic Books in the Fall of 2013. \\ For more
  information about the film, visit http://ritesofloveandmath.com}

\author{Edward Frenkel}

\address{Department of Mathematics, University of California,
Berkeley, CA 94720, USA}

\maketitle


The lights were dimmed... After a few long seconds of silence the
movie theater went dark. Then the giant screen lit up, and black
letters appeared on the white background:

\bigskip

\begin{center}

Red Fave Productions

\bigskip

in association with Sycomore Films

\bigskip

with support of

\bigskip

Fondation Sciences Math\'ematiques de Paris

\bigskip

present

\bigskip

{\bf Rites of Love and Math}

\end{center}

\bigskip

The 400-strong capacity crowd was watching intently. I'd seen it
countless times in the editing studio, on my computer, on TV... But
watching it for the first time on a panoramic screen was a special
moment which brought up memories from the year before.

\bigskip

I was in Paris as the recipient of the first {\em Chaire d'Excellence}
awarded by Fondation Sciences Math\'ematiques de Paris, invited to
spend a year in Paris doing research and lecturing about it.

Paris is one of the world's centers of mathematics, but also a
capital of cinema. Being there, I felt inspired to make a movie about
math. In popular films, mathematicians are usually portrayed as
weirdos and social misfits on the verge of mental illness, reinforcing
the stereotype of mathematics as a boring and irrelevant subject, far
removed from reality. Would young people want a career in math or
science after watching these movies? I thought something had to be
done to confront this stereotype.

My friend, mathematician Pierre Schapira, introduced me to a young
talented film director, Reine Graves. A former fashion model, she had
previously directed several original, bold short films (one of which
won the Pasolini Prize at the Festival of Censored Films in Paris). At
a lunch meeting arranged by Pierre, she and I hit it off right away. I
suggested we work together on a film about math, and she liked the
idea. Months later, asked about this, she said that she felt
mathematics was one of the last remaining areas where there was
genuine passion.

As we started to brainstorm what our film would be about, I showed
Reine a couple of photographs I had made, in which I painted tattoos
of mathematical formulas on human bodies. We decided we would try
to make a film involving the tattoo of a formula.

Tattoo, as an art form, originated in Japan. I had visited Japan a
dozen times, was fascinated by the Japanese culture.  We turned to the
Japanese cinema for inspiration; in particular, to a film by the great
Japanese writer Yukio Mishima {\em Rite of Love and Death}, based on
his short story {\em Y\^ukoku} (or {\em Patriotism}). Mishima himself
directed and starred in it.

This film made a profound impression on me. It was as though I was
possessed by a powerful force.

{\em Rite of Love and Death} is black-and-white; it unfolds on the
austere stylized stage of the Japanese Noh Theater. No dialogue, with
music from Wagner's opera {\em Tristan and Isolde} playing in the
background. There are two characters: a young officer of the Imperial
Guard, Lieutenant Takeyama, and his wife, Reiko. The officer's friends
stage an unsuccessful {\em coup d'etat} (here the film refers to
actual events of February 1936, which Mishima thought had a dramatic
effect on Japanese history). The Lieutenant is given the order to
execute the perpetrators of the {\em coup}, which he cannot do -- they
are close friends. But neither can he disobey the order of the
Emperor. The only way out is ritual suicide, {\em seppuku} (or {\em
  harakiri}). When he tells Reiko, she says she will follow him to the
better world. After they make love, the Lieutenant commits {\em
  seppuku} (this is shown in graphic detail). Then Reiko kills herself
by driving a knife in her throat... At the end of the film we see them
both lying dead in final embrace on the beautifully groomed pebbles of
a traditional Zen garden.

The 29 minutes of film touched me deeply. I could sense the vigor and
clarity of Mishima's vision. His presentation was forceful, raw,
unapologetic. You may disagree with his ideas (and in fact his vision
of the intimate link between love and death does not appeal to
me), but you have to respect the author for being so strong and
uncompromising.

Mishima's film went against the usual conventions of cinema: it was
silent, with written text between the ``chapters'' of the movie to
explain what's going to happen next. It was theatrical; scenes
carefully staged, with little movement. But I was captivated by the
undercurrent of emotion underneath. (I did not know yet the details of
Mishima's own death, its eerie resemblance to what happened in his
film -- and this was probably for the best.)

Perhaps, the film resonated with me so much in part because Reine and
I were also trying to create an unconventional film, to talk about
mathematics the way no one had talked about it before. I felt that
Mishima had created the aesthetic framework and language that we were
looking for.

I recount what happened next in my forthcoming book {\em Love and
  Math}. I called Reine and told her that we should make a film just
like Mishima's. ``But what will our film be about?'' she asked.
Suddenly, words started coming out of my mouth. Everything was crystal
clear.

\medskip

``A mathematician creates a formula of love,'' I said, ``but then
discovers the flip side of the formula: it can be used for Evil as
well as for Good. He realizes he has to hide the formula to protect it
from falling into the wrong hands. And he decides to tattoo it on the
body of the woman he loves.''

\medskip

We decided to call our film {\em Rites of Love and Math}. We
envisioned it as an allegory, showing that a mathematical formula can
be beautiful like a poem, a painting, or a piece of music. The idea
was to appeal not to the cerebral, but to the intuitive, visceral. Let
the viewers {\em feel} rather than {\em understand} it first. We
thought that this would make mathematics more human, inspire viewer's
curiosity about it.

We also wanted to show the passion involved in mathematical
research. People tend to think of mathematicians as sterile, cold. But
the truth is that our work is full of passion and emotion. And the
formulas you discover really do get {\em under your skin} -- that was
the intended meaning of the tattooing in the film.

\medskip

In the film, the Mathematician discovers the ``formula of love.'' Of
course, this is a metaphor: We are always trying to reach for complete
understanding, ultimate clarity, want to know everything. In the real
world, we have to settle for partial knowledge and understanding. But
what if someone were able to find the ultimate Truth; what if it could
be expressed by a mathematical formula?  This would be the ``formula
of love.''

Such a formula, being so powerful, must also have a flip side: it
could also be used for evil. This is a reference to the dangers of
modern science. Think of a group of theoretical physicists trying to
understand the structure of the atom. What they thought was pure
scientific research inadvertently led them to the discovery of atomic
energy. It brought us a lot of good, but also destruction and
death. Likewise, a mathematical formula discovered as part of our
quest for knowledge could potentially lead to disastrous consequences.

So our protagonist, the Mathematician, hides his formula by tattooing
it on the body of the woman he loves. It's his gift of love, the
product of his creation, passion, imagination. But who is she? In the
framework of the mythical world we envisioned, she is the incarnation
of Mathematics, Truth itself (hence her name Mariko, ``truth'' in
Japanese; and that's why the word {\em ISTINA}, ``truth'' in Russian,
is calligraphed on the painting hanging on the wall). The
Mathematician's love for her is meant to represent his love for
Mathematics and Truth, for which he sacrifices himself. But she
survives and carries his formula, as she would their child. The Truth
is eternal.

\medskip

As Reine and I were getting more excited about our project, so did
people around us. Soon, a crew of about 30 people was working on the
film. Rapha\"el Fernandez, a talented musician, composed original
music. We ordered a kimono and a painting. An artist was working on
the decor. The film was taking life of its own.

\medskip

The shooting took three days, in July of 2009. Those were some of the
most exciting, and exhausting, days of my life. I wore several hats:
co-director (with Reine), producer, actor... All of this was new to
me, and I was learning on the job. It was an amazing journey, a
wonderful collaboration with Reine and other filmmakers and artists
helping us to fulfill our dream.

\medskip

The central scene of the movie is the making of the tattoo. We shot it
on the last day. Since I never had a tattoo, I had to learn about the
process. These days tattoos are made with a machine, which would be an
anachronism on the stage of Noh theater. In the past, tattoos were
engraved with a bamboo stick -- a longer, more painful process. I've
been told it's still possible to find tattoo parlors in Japan which
use this old technique. This is how we presented it in the film.

\bigskip

\begin{center}

\includegraphics[height=75mm]{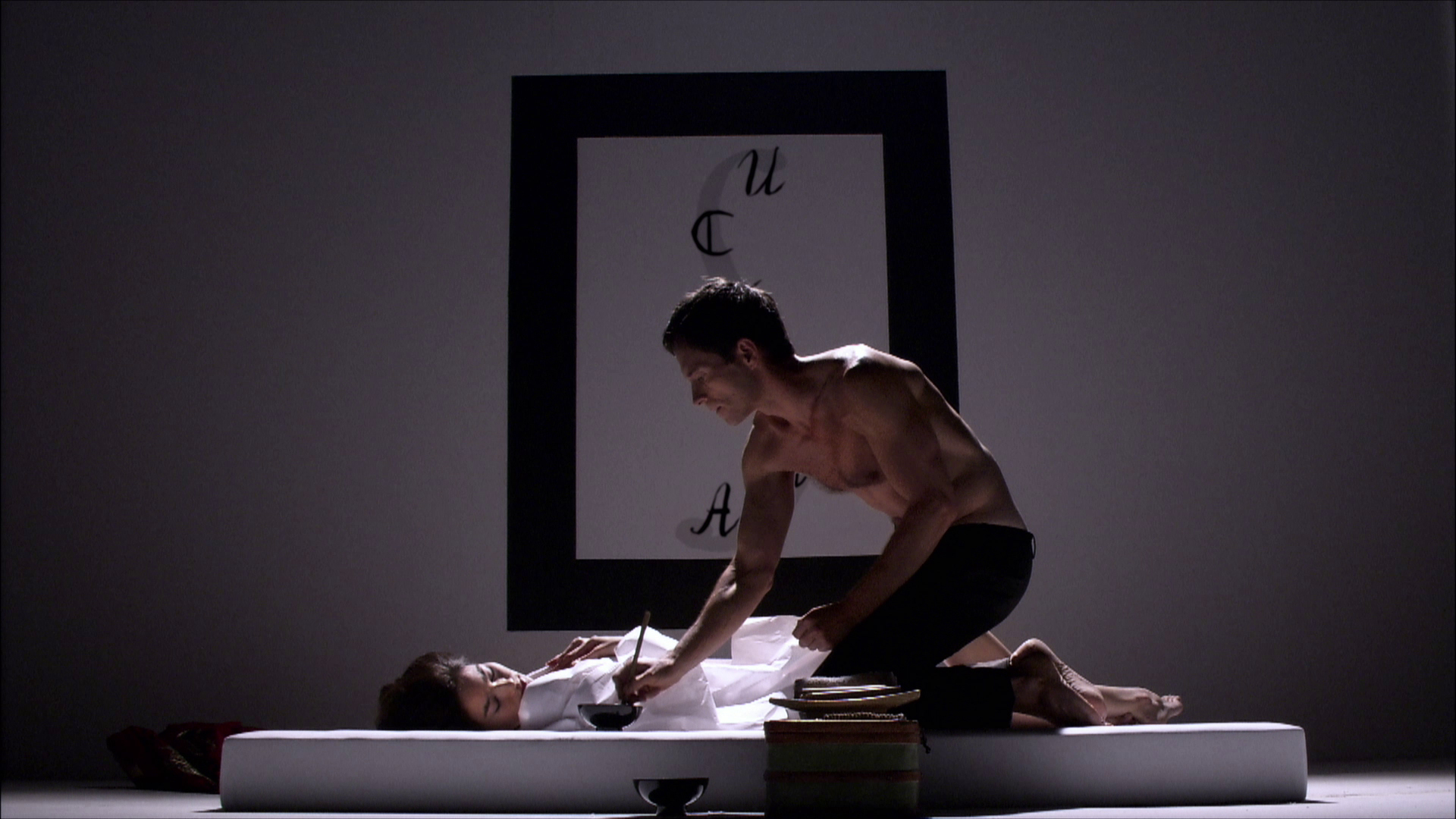}

\end{center}

\bigskip

Oriane Giraud, our special effects artist, asked me a few days before
the shooting to give her the formula that would be tattooed, so she
could create the blueprint. Which formula should play the role of
``formula of love''? A big question! It had to be sufficiently
complicated (it's a formula of love, after all), aesthetically
pleasing. We wanted to convey that a mathematical formula could be
beautiful in content as well as form. And I wanted it to be {\em my}
formula.

Doing ``casting'' for the formula of love, I stumbled on this:

\begin{align*}    \notag
\int_{\pone} \omega F(qz,\ol{q}\zb) &= \sum_{m,\ol{m} = 0}^\infty
\underset{|z|<\ep^{-1}}\pf \omega_{z\zb} \ z^{m} {\zb}^{\ol{m}} dz
  d\zb \; \cdot \left. \frac{q^m \ol{q}^{\ol{m}}}{m! {\ol{m}}!}
{\p}_z^{m} {\p}_{\zb}^{\ol{m}} F \right|_{z=0} \\ &+
q \ol{q} \sum_{m,\ol{m} = 0}^{\infty} \left. \frac{q^m
\ol{q}^{\ol{m}}}{m!  {\ol{m}}!} {\p}_w^{m} {\p}_{\wb}^{\ol{m}}
\omega_{w\wb} \right|_{w=0} \cdot \underset{|w|<q^{-1}\ep^{-1}}\pf F \
w^{m} {\wb}^{\ol{m}} dw  d\wb.
\end{align*}

It appears as formula (5.7) in a 100-page paper \cite{FLN}, {\em
  Instantons Beyond Topological Theory I}, which I wrote with my two
good friends, Andrey Losev and Nikita Nekrasov, in 2006.\footnote{In
  the version published in {\em Journal de l'Institut de
    Math\'ematiques de Jussieu}, there is a footnote explaining that
  it played the role of ``formula of love'' in {\em Rites of Love and
    Math}.}

When we show the film, people always ask: What does this formula mean?
Which is exactly what we were hoping for. If we had made a film in
which I wrote this formula on a backboard and tried to explain its
meaning, how many people would care for it? But seeing it in the form
of a tattoo elicited a totally different reaction. It really got under
everyone's skin.

\bigskip

So {\em what does it mean}? This work was the first installment in a
series of papers we wrote about a new approach to quantum field
theories with ``instantons'' -- special solutions of the theory
minimizing the ``action.'' Quantum field theory is a mathematical
formalism for describing the behavior and interaction of elementary
particles. Though it has been successful in accurately predicting a
wide range of phenomena, there are many fundamental issues that are
still poorly understood. Let me recall that all atoms consists of
protons, neutrons, and electrons. Protons and neutrons, in turn,
consist of smaller particles, called quarks. And those quarks are
confined there -- they cannot be separated. A proper theoretical
explanation of this phenomenon is still lacking.

In the conventional (so-called, perturbative) approach to quantum
field theory the starting point is the so-called ``free theory''
describing idealized non-interacting particles. Then we ``turn on''
the interaction between them. The problem is that in this approach the
contribution of each instanton appears to be negligible, even though
altogether they may ``conspire'' to create a powerful effect. Many
physicists believe that the difficulties of the conventional formalism
with taking the instantons into account could be the reason we don't
have a satisfactory explanation of the confinement of quarks and other
similar effects.

In our paper, we proposed a new approach, in which the starting point
is not a free theory, but an idealized {\em interacting} theory in
which the instantons are present from the beginning. The advantage of
out theory is that the main quantities -- the so-called ``correlations
functions'' -- are expressed by {\em finite-dimensional} integrals (in
contrast to the conventional formalism, in which the integrals are
infinite-dimensional and hence poorly defined). Therefore, our theory
is in principle completely solvable.

The above formula expresses the identity between two ways to compute
correlation functions in our theory. We discovered it when we were
working on this project three years earlier, in April of 2006, also in
Paris.

Our world is four-dimensional (if we include both space and time
dimensions), but four-dimensional theories are very complex. To
simplify matters, we looked at the analogous two-dimensional, and then
one-dimensional, models. That is to say, there is only time
dimension. Such models, commonly referred to as ``quantum mechanics,''
describe a single particle moving in a particular space (which could
be of any dimension). Despite the simplifications, these
one-dimensional models possess the salient features of the more
realistic, four-dimensional, models. That's why it is useful to study
them.

We considered the theory in which a particle was moving on the sphere,
also known in mathematics as the ``complex projective line,'' denoted
$\pone$ (${\mathbb P}$ stands for ``projective,'' and $1$ for
one-dimensional, as a complex manifold). You can see this notation
under the integral sign on the left hand side of the formula.

We wanted to compute the simplest correlation function in this theory,
involving two ``observables,'' denoted by $F$ and $\omega$ in our
formula. On the one hand, the answer is given by an integral over the
sphere; this is the left hand side of the formula. On the other hand,
in our new theory we were getting a different answer: a sum over the
``intermediate states,'' appearing on the right hand side. This answer
is surprising, and so is the equality between the two expressions. If
our new approach were correct, the two sides would have to be equal to
each other. And indeed they are.

In other words, our formula says that two ways of computing the
correlation functions -- the old and the new -- give the same
answer. Little did we know at the time we discovered it that it would
soon be slated to play the role of formula of love.

\bigskip

Oriane liked the formula, but said it was too involved for a
tattoo. Could I simplify it? I slighly changed the notation, and
here's how it appears in our film:

\bigskip

\begin{center}

\includegraphics[height=75mm]{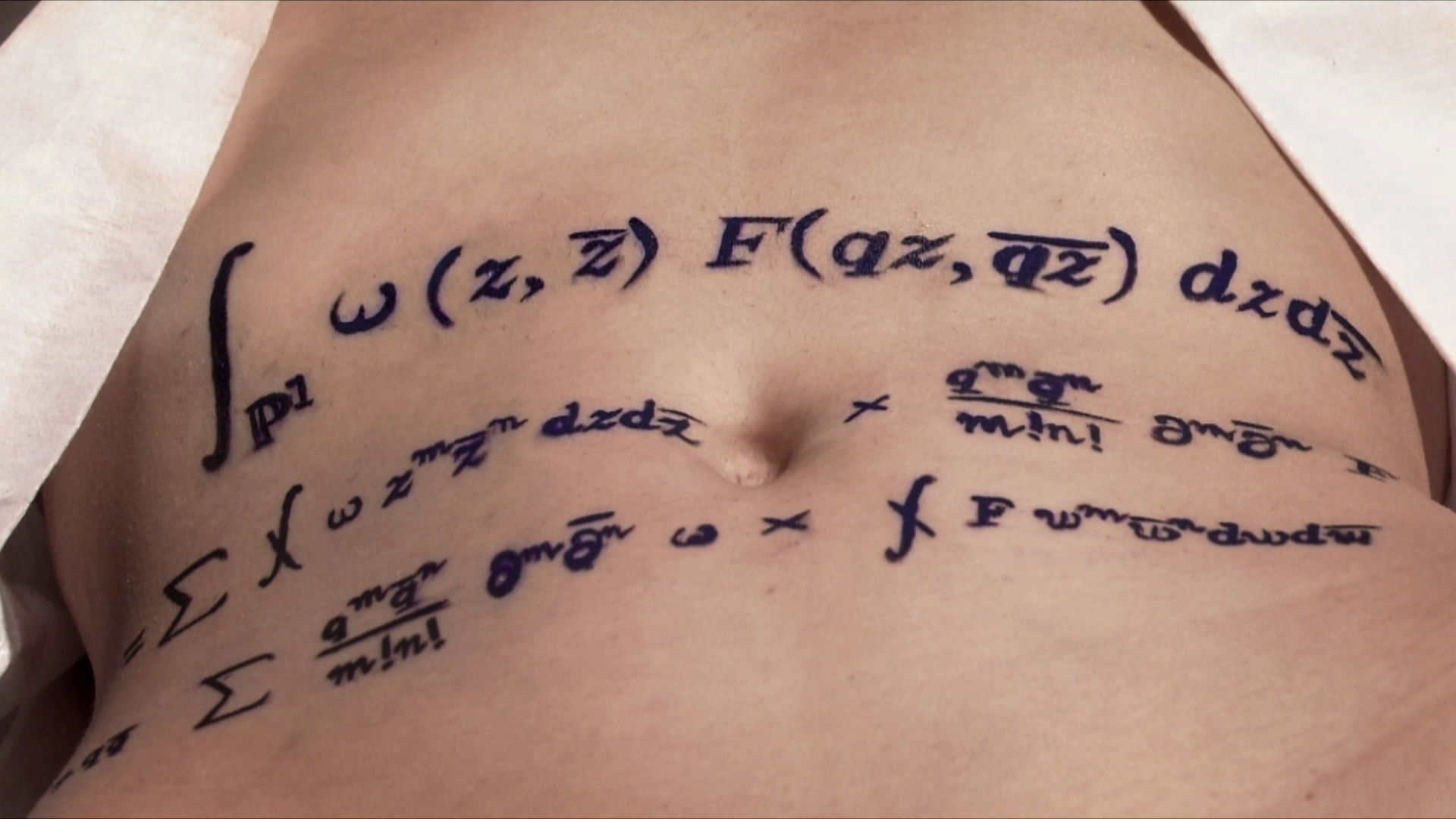}

\end{center}

\bigskip

The work on the tattoo scene took us many hours. It was
psychologically and physically draining both for me and for Kayshonne
Insixieng May, the actress playing Mariko. We finished shooting close
to midnight. It was an emotional moment for all of us on the set,
after everything we had been through together.

\medskip

Two months of post-production followed, at {\em Sycomore Films}, a
French produc\-tion company, with our multi-talented editor
Thomas Bertay and visual effects magician Pierre Borde. And finally,
the film finished, it was time to organize the premiere. Folks from
Fondation Sciences Math\'ematiques de Paris, which generously
supported our film, agreed to sponsor the premiere. We picked a
wonderful venue: {\em Max Linder Panorama}. An old movie theater with
a huge screen and modern video and sound systems, it is one of best
theaters in Paris.

The crowd that gathered at the theater on April 14, 2010 was diverse:
mathematicians, artists, filmmakers... Most of our crew was there,
including Reine, Kayshonne, and our Director of Photography {\em
  extraordinaire} Daniel Barrau. On the big screen the picture was
sharp and crisp, the colors vivid (the decision to use the most
expensive Sony camera had paid off).

\medskip

The first articles about the film started to appear. {\em Le Monde}
called {\em Rites of Love and Math} ``a stunning short film'' that
``offers an unusual romantic vision of mathematicians.'' And the {\em
  New Scientist} wrote:

\medskip

{\em It is beautiful to look at... If Frenkel's goal was to bring more
  people to maths, he can congratulate himself on a job well done. The
  formula of love, which is actually a simplified version of an
  equation he published in a 2006 paper on quantum field theory
  entitled ``Instantons beyond topological theory I'', will probably
  soon have been seen -- if not understood -- by a far larger audience
  than it would otherwise ever have reached.}

\bigskip

Since then, the film has been shown at international film festivals in
France, Spain, and in Berkeley, California. There have been more
showings in Paris, Kyoto, Madrid, Santa Barbara, Bilbao, Venice... The
screenings and the ensuing publicity gave me the opportunity to meet
many people and hear different opinions. At first, this came as
culture shock. Some of my mathematical works can be fully understood
only by a small number of people; sometimes, no more than a dozen in
the whole world at first. This was a film intended for a wide
audience: hundreds were being exposed to it. And of course, they all
interpreted it in their own ways.

In mathematics there is only one truth, and only one path to reach
that truth. My mathematical work is perceived and interpreted in
essentially the same way by everybody who reads it. Not so in cinema
and in the arts in general. First, there isn't a single truth, and
second, there are so many different paths to express the truth. And
the viewer is always part of an artistic project: at the end of the
day, it's all in the eye of the beholder. I have no influence over
their perception. Coming to terms with this was a challenge to me, but
gradually I came to embrace it. We get enriched when others share
their views and insights. What matters the most for a work of art is
that it touches people in some way, does not leave them indifferent.

\bigskip

When we show the film, people invariably ask: ``Do you know the
formula of love?''  My response: ``Every formula we create is a
formula of love.''  Doing mathematics is a creative pursuit that
requires passion, just like painting, music, and poetry. In order to
discover something new and eternal about the world, you have to be in
love with what you do.

\newpage

\begin{center}

\includegraphics[height=190mm]{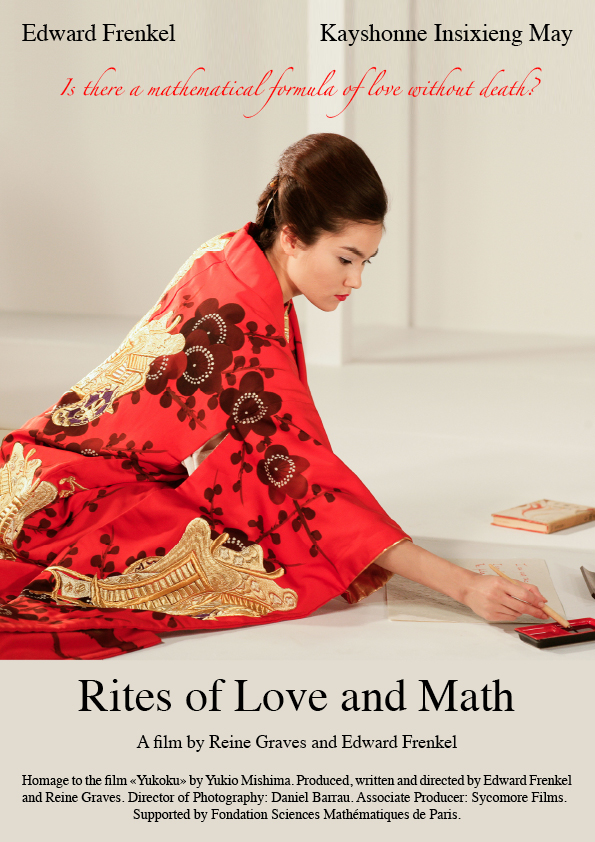}

\bigskip

{\small Film poster}

\end{center}

\bigskip

\section*{Film Credits}

{\em Rites of Love and Math}, 2010. 26 minutes, in color.

Written, produced, and directed by Edward Frenkel and Reine Graves.

Associate Producer {\em Sycomore Films}.

Cinematography by Daniel Barrau.

Music by Richard Wagner and Rapha\"el Fernandez.

With Edward Frenkel and Kayshonne Insixieng May.

More information: http://ritesofloveandmath.com

\bigskip

\bigskip

{\em Rite of Love and Death} ({\em Y\^ukoku}), 1965. 29 minutes, in
black-and-white.

Written, produced, and directed by Yukio Mishima.

Associate producer Haraoki Fujii.

Associate director Masaki Domoto.

Cinematography by Kimio Watanabe.

Music by Richard Wagner.

With Yukio Mishima and Yushiko Tsuruoka.

\vspace*{10mm}

\section*{Acknowledgments}

I thank Michele Emmer for the invitation to show our film at the
annual Symposium ``Matematica e Cultura'' in Venice and to contribute
to the Symposium Proceedings. I am grateful to Thomas Farber for his
comments on a draft of this article.

\end{document}